\documentclass[12pt]{amsart}
\usepackage{amssymb,latexsym, amsmath, amsxtra}
\usepackage{graphicx}
\newtheorem{conjecture}{Conjecture}
\newtheorem{theorem}{Theorem}

\begin{document}

\title{In support of  $n$-correlation}
 
\abstract 
In this paper we examine  
$n$-correlation for either the eigenvalues of a unitary group of
random matrices or for the zeros of a unitary family of $L$-functions
in the important situation when the correlations are detected via
test functions whose Fourier transforms have limited support. This
problem first came to light in the work of Rudnick and Sarnak
in their study of the $n$-correlation of zeros of a fairly general
automorphic $L$-function. They solved the simplest instance of this
problem when the total support was most severely limited, but had to
work extremely hard to show their result matched random matrix theory
in the appropriate limit. This is because they were comparing their
result to the familiar determinantal expressions for $n$-correlation
that arise naturally in random matrix theory.  In this paper we deal
with arbitrary support and show that there is another expression for
the $n$-correlation of eigenvalues that translates easily into the
number theory case and allows for immediate identification of 
which terms survive the restrictions placed on the support of the
test function.

 \endabstract

\author {J.B. Conrey}
\address{American Institute of Mathematics,
360 Portage Ave, Palo Alto, CA 94306 USA and School of
Mathematics, University of Bristol, Bristol, BS8 1TW, United
Kingdom} \email{conrey@aimath.org}

\author{N.C. Snaith}
\address{School of Mathematics,
University of Bristol, Bristol, BS8 1TW, United Kingdom}
\email{N.C.Snaith@bris.ac.uk}

\subjclass[2010]{Primary 11M26; Secondary 11M50}

\keywords{$n$-correlation, Riemann zeta-function, random matrix theory}

\thanks{
Research of the first author supported by the American Institute
of Mathematics and by a grant from the National Science
Foundation. The second author was sponsored by the Air
Force Office of Scientific Research, Air
  Force Material Command, USAF, under grant number FA8655-10-1-3088.} \maketitle

\maketitle

\section{Introduction}
In 1973 Hugh Montgomery published his revolutionary work [Mon] on the zeros of the Riemann zeta-function.
Assuming the Riemann Hypothesis, so that the zeros are  $1/2+i \gamma$, he evaluated asymptotically sums which are basically of the form
$$\sum_{0 < \gamma \le \mathcal{T}\atop 0\le \gamma'\le \mathcal{T}} f((\gamma-\gamma')(\log \mathcal{T})/2 \pi)$$
where $f$ is a function whose Fourier transform is supported on $(-1,1)$. 
His analysis led him to make the   
\begin{conjecture}
For fixed $\alpha < \beta$,
\begin{eqnarray}
\sum_{{0< \gamma\le \mathcal{T}\atop 
0< \gamma'\le \mathcal{T}}
\atop
2\pi \alpha /\log \mathcal{T} \le \gamma-\gamma'\le 2 \pi \beta /\log \mathcal{T}}
1 \sim \int_\alpha^\beta\left(1-\left(\frac{\sin \pi u }{\pi u}\right)^2~du +\delta(\alpha, \beta)\right) \frac{\mathcal{T}}{2\pi} \log \mathcal{T}
\end{eqnarray}
as $\mathcal{T}$ tends to infinity. Here $\delta(\alpha,\beta)= 1$ if $0\in [\alpha,\beta]$, $\delta(\alpha, \beta)=0$ otherwise. 
\end{conjecture}
Montgomery credits Freeman Dyson with the observation that this conjecture is consistent with the zeta-zeros being spaced like the eigenvalues of 
large random Hermitian or unitary matrices. 
Montgomery goes on to say,
\begin{quote}
 One might extend the present work to investigate the k-tuple correlation of the zeros of the zeta function. 
If the analogy with random complex Hermitian matrices appears to continue, then one might conjecture that the $k$-tuple correlation 
function $\hat{F}(u_1, u_2, ..., u_k)$ is given by
\begin{equation*}
\hat{F}(u_1, u_2, ..., u_k)=\det A\end{equation*}
where $A = [a_{ij}]$ is the $k \times k$ matrix with entries $a_{ii} = 1, a_{ij} = (\sin \pi (u_i -u_j))/\pi(u_i -u_j)$ for  $i\ne j$.
\end{quote}

We can explain rather simply the random matrix theorem behind all of this. 
Let $U(N)$ denote the group of $N\times N$ unitary matrices and let $dU$ signify the Haar measure for this group. 
For a unitary matrix $U$ we let $\theta_1,\dots , \theta_N$ denote the $N$ eigenangles, i.e $e^{i\theta_j}$ are the eigenvalues. 
A beautiful theorem (see [Meh]) in random matrix theory gives the $n$-correlation of these eigenvalues:
\begin{theorem} \label{thm:meh}
Let $F:[0,2\pi]^n\to \mathbb C$ be an integrable function of
$n$-variables. Then
\begin{eqnarray*}
\int_{U(N)}\sideset{}{^*}\sum_{(j_1,\dots , j_n)}
F(\theta_{j_1},\dots ,\theta_{j_n}) dU = \frac{1}{(2\pi)^n}
\int_{[0,2\pi]^n} F(\alpha_1,\dots,\alpha_n)
  \det_{n\times n} S_N(\alpha_k-\alpha_j)~d\alpha_1\dots
  ~d\alpha_n,
\end{eqnarray*}
where $\sideset{}{^*}\sum$ indicates that the sum is over $n$-tuples of distinct indices $(j_1,\dots,j_n)$ with $1\le j_i\le N$, i.e. 
$j_i\ne j_k$ if $i\ne k$,  and where 
\begin{eqnarray}
S_N(\alpha)=\frac{\sin \frac {N\alpha}{2}}{\sin \frac \alpha 2 }.
\end{eqnarray}
\end{theorem}

Indeed, Rudnick and Sarnak [RuSa, 1996] did take up the question of the $n$-correlation of the zeros of the Riemann zeta-function,
and so studied sums basically of the form
\begin{eqnarray}
\label{eqn:RSsum}
\sum f(\gamma_1 (\log \mathcal{T}) /2 \pi,\dots, \gamma_n (\log \mathcal{T})/2 \pi)
\end{eqnarray}
where again the support of the Fourier transform of $f$ is an issue; the proof requires that
$f(u_1,\ldots,u_n)$ be supported in $\sum|u_j|<2$.
 (In fact Rudnick and Sarnak study these sums for an {\it arbitrary} $L$-function.)
But there is a surprising turn of events in [RuSa]; namely at the end of the natural series of calculations, it is not clear whether 
or not the answer they get is consistent with the random matrix theory prediction!
To resolve this issue requires a rather large detour into some combinatorics which seem to be unrelated to their central pursuit.

In the calculation by Rudnick and Sarnak [RuSa] of the $n$-correlation of the 
zeros of the Riemann zeta-function,
the inherent difficulty in evaluating the sum (\ref{eqn:RSsum}) occurs
 in 
 dealing with off-diagonal terms in the prime number sums that arise after the use of the {\it explicit formula}.
The limitations in the support of the Fourier transform of the test function $f$ allow one to ignore such considerations.
  These off-diagonal terms are the subject of the work of  Bogomolny and Keating [BoKe] who  construct  a beautiful demonstration
 that the combinatorics of prime number sums 
arising from the Hardy-Littlewood prime pair conjectures lead to the determinantal formula of random matrix theory.

As mentioned, in [RuSa], after deriving an expression for the $n$-correlation in the situation that their test function
has a Fourier transform with   support limited to a proper subset of $(-2,2)$\footnote{
In the case of pair-correlation, the  test functions in [RuSa] are written with two variables whereas in [Mon] one variable is used. 
This is the source of the apparent discrepancy in ranges of support of the Fourier transforms where 
[Mon] has $(-1,1)$ and [RuSa] has $(-2,2)$; the results are of the same quality.}
, the authors are faced with the non-trivial task of verifying that
their answer agrees with theorems from random matrix theory. 
A similar problem occurs in two works on the $n$-level density of zeros of quadratic
$L$-functions, first considered  by Mike Rubinstein [Rub], and subsequently by Peng Gao in his thesis.
In the former case Rubinstein finds, as do Rudnick and Sarnak,  an adhoc method to verify the  consistency
of the number theory and the random matrix theory calculations.
Peng Gao is unable to make this verification in his situation. A subsequent paper [ERR] 
rectifies this situation and completes the verification in Gao's case by a very clever  appeal to 
 zeta-functions over function fields. 

In light of this discussion it is natural to ask about
\begin{eqnarray}&& 
\int_{U(N)}\sideset{}{^*}\sum_{(j_1,\dots , j_n)}
f\left(\frac{N\theta_{j_1}}{2\pi},\dots ,\frac{N\theta_{j_n}}{2\pi}\right) dU 
\end{eqnarray}
and the scaled limit of this quantity as $N\to \infty$ in the situation where $f$ is a  translation invariant function that has a 
Fourier transform with limited support. This is what we will study in this paper. 
  
The usual proof of Theorem \ref{thm:meh} is very short and elegant; it relies on  what is known as Gaudin's lemma. 
In our previous article [CoSn]  we give a different derivation  that is far more complicated, but which has the one advantage
that it can be copied step by step to 
obtain an analogous formulation for the $n$-correlation of zeros of the $\zeta$-function complete 
with all of the arithmetic terms which, of course, do not appear in random matrix theory. 
Our starting point for the alternate proof of the random matrix theorem above is 
the unitary Ratios Theorem, see [CFS], and our starting point for the conjectural detailed zeta-function analogue
are a collection of conjectures known as the  $L$-functions Ratios Conjectures; see [CoSn1] 
for an introduction to this useful set of conjectures.   

In this paper, we show that the formula for $n$-correlation derived in [CoSn] allows for immediate identification of the terms which survive a restriction on the support of the test function, and so makes the verification of Rudnick and Sarnak's formula in the unitary case straightforward, 
eg. for the $n$-correlation
of the zeros of the Riemann zeta-function. (In [RuSa] more general $L$-functions are considered; here we focus on the case of 
the Riemann zeta-function, but higher degree $L$-functions don't present any extra difficulty.)  Our proof is natural in that  it 
explains the situation for a test function whose Fourier transform has any range of support. 
Indeed, the work of [CLLR] uses the asymptotic large sieve to investigate the pair correlation of the zeros of all primitive Dirichlet 
$L$-functions; this is a unitary family and the allowable test functions can have double the range of support as in the Rudnick-Sarnak 
situation. If these co-authors extend their work to $n$-correlation then their situation compared to Rudnick-Sarnak's will be exactly as Peng 
Gao's situation compared to Rubinstein's, i.e. virtually intractable though ordinary combinatorial ideas. However, when this 
happens, the crucial step needed to conclude that their answer agrees with RMT will be furnished by this paper. 

In her thesis, Amy Mason has derived the
$n$-correlation functions for the orthogonal and symplectic groups in a work analogous to
 [CoSn] which is about the unitary group. The techniques of the  
the current paper may well extend to the situation where one has orthogonal or symplectic symmetry. 
If so, that work will give an alternate proof of the result of [ERR] that Peng Gao's theorem matches 
random matrix theory.


\section{Statement of results}
For simplicity we state the results of Rudnick and Sarnak  for the $n$-correlation of the zeros $1/2+i\gamma_j$ of the Riemann
zeta-function, though there is no difficulty working in their full generality. 
Before doing so, we introduce their vector notation  as a convenient way to express the combinatorial sum
that arises in their work. 
Let
\begin{eqnarray}
\left\{
\begin{array}{ll}
\mathbf e_{i,j} &= \mathbf e_i -\mathbf e_j\\
\mathbf e_i & = (0,\dots,1, \dots ,0) \mbox{ the $i$th standard basis vector} 
\end{array}
\right.
\end{eqnarray}

\begin{theorem}[]  [RuSa] Theorem 3.1. Let $h_j$, for $1\le j\le n$, be  rapidly decaying functions with 
\begin{eqnarray} \label{eqn:gtoh}
h_j(x)=\int_{\mathbf R} g_j(t)e^{ixt} ~ dt
\end{eqnarray} 
where $g_j$ is smooth and compactly supported. Let $e(z)=e^{2\pi i z}$ and let $\delta$ be the Dirac $\delta$-function.  
Suppose that
\begin{eqnarray} \label{eq:C1}
f(x_1,\dots , x_n)=\int_{{\mathbf R}^n}\Phi(\xi_1,\dots \xi_n) \delta(\xi_1+\dots +\xi_n)e(-x_1\xi_1-\dots -x_n\xi_n) d\xi_1 \dots d\xi_n
\end{eqnarray}
where $\Phi$ is smooth, even, and compactly supported in such a way that 
\begin{eqnarray} \label{eq:C2}
\Phi(\xi_1,\dots ,\xi_n)=0
\end{eqnarray}
whenever
\begin{eqnarray} \label{eq:C3}
|\xi_1|+\dots+|\xi_n| >2-\epsilon 
\end{eqnarray}
for some fixed $\epsilon >0$.  Let $\mathcal L = \log \mathcal{T}$.
 Then
\begin{eqnarray*} &&
 \sum_{\gamma_1,\dots ,\gamma_n}h_1\left(\frac{\gamma_1}{\mathcal{T}}\right) \dots
h_n\left(\frac{\gamma_n}{\mathcal{T}}\right)f\left(\frac{\mathcal L\gamma_1}{2\pi},\dots, \frac{\mathcal L \gamma_n}{2\pi}
\right)\\
&& \quad  =  \kappa(\mathbf h)\frac{\mathcal{T}\mathcal L}{2\pi} \bigg( \Phi(0)  +\sum_{r=1}^{[n/2]}\sum_{i(t)<j(t)\atop t\le r}
\int |v_1|\dots |v_r| \Phi(v_1 \mathbf e_{i(1),j(1)}+\dots + v_r \mathbf e_{i(r),j(r)}) ~ dv 
\bigg) +O(\mathcal{T}),
\end{eqnarray*}
where   the sum is over all disjoint pairs of indices $i(t) < j(t)$ in $\{ 1,2 ,\dots , n\}$ when $t\le r$ 
and where 
\begin{eqnarray}
\kappa(\mathbf h) =\int_{\mathbf R} h_1(u)\dots h_n(u) ~du.
\end{eqnarray}
\end{theorem}

The main result of this paper is the random matrix  analogue of this theorem. 
For the purposes of making a close connection with number theory it is convenient to let the eigenangles ``wrap around.''
 Thus, for an $N \times N$ 
 unitary matrix
 $X$ we will let its
eigenangles be
\begin{eqnarray}
\dots \le \theta_{-R} \le  \theta_{-R+1}\le \dots \le  \theta_0 \le \theta_1 \le \dots \le \theta_R \le  \dots
\end{eqnarray}
where
\begin{eqnarray}
\theta_{r+kN}=\theta_r+ 2\pi k.
\end{eqnarray}  
Now we have infinitely many eigenangles just as for the Riemann zeta-function we have infinitely many zeros. 
We are interested in the sum
\begin{eqnarray}
\sum_{j_1,\dots,j_n}h_1\left(\frac{\theta_{j_1}}{\mathcal{T}}\right)\dots
 h_n\left(\frac{\theta_{j_n}}{\mathcal{T}}\right) f\left(\frac{N\theta_{j_1}}{2\pi},\dots ,\frac{N\theta_{j_n}}{2\pi}\right)
\end{eqnarray}
 where each $j_k$ now runs over all integers.  
 This sum exactly parallels the above sum over zeta-zeros. 

We will assume throughout this paper that $n$ is fixed and that $N$ is sufficiently
large in terms of $n$. Also, we regard $\mathcal T$ as sufficiently large in terms of $N$.
 
\begin{theorem}  \label{th:r1}   With the same conditions on $h$ and $f$ (and hence on $\Phi$), we have 
\begin{eqnarray} &&
  \int_{U(N)}\sum_{j_1,\dots ,j_n}h_1\left(\frac{\theta_{j_1}}{\mathcal{T}}\right)\dots
 h_n\left(\frac{\theta_{j_n}}{\mathcal{T}}\right)
f\left(\frac{N \theta_{j_1} }{2\pi}, \dots ,\frac{N \theta_{j_n}}{2\pi}\right)
 dU\nonumber \\
&&\qquad =\kappa(\mathbf h)\frac{N\mathcal{T}}{2\pi}  \sum_{K+L+M=
\{1,\dots,n\}\atop |K|=|L|} \sum_{\sigma\in S_{|K|}} 
\int_{  \xi_{k_j}>0 }
\xi_{k_1}\dots \xi_{k_K}
 \Phi\left(\sum_{j=1}^K \xi_{k_j}\mathbf e_{k_j,\ell_{\sigma(j)}}\right) d\xi +O(N).
\end{eqnarray}
The sum $K+L+M$ is over set partitions of $\{1,2,\dots,n\}$ into 3 parts, but empty sets are allowed.
Also, $S_{|K|}$ is the symmetric group of order $|K|!$.
\end{theorem}

This theorem is slightly more general than the theorem above of Rudnick and Sarnak and has a smaller error term;
but the main terms exactly match when $\Phi$ is assumed to be even as in [RuSa] and $N$ takes the place of $\log \mathcal{T}$.
As a consequence, we have proven 
   that the result of Rudnick and Sarnak 
 agrees with random matrix theory without going through the 
complex combinatorial considerations that they undergo in the second half of their paper.

 Let $q>0$. 
In general, one could ask to compare  number theory and random matrix theory 
to see what happens if one works with a set of test functions $\mathcal T_q$ which are defined as above except
that the support of $\Phi$ is restricted to
\begin{eqnarray}
|\xi_1|+\cdots + |\xi_n|<2 q.
\end{eqnarray}

We have developed an approach in [CoSn] that treats number theory  and random matrix theory in parallel.  From this work it is a simple 
matter to read off the results when dealing with test functions from $\mathcal T_q$. This result is described below.

\section{Eigenvalue correlations}
 
Before stating the theorem of [CoSn] we describe some notation.
We let
\begin{eqnarray}\label{eq:z}
z(x)=\frac{1}{1-e^{-x}}.
\end{eqnarray}
In our formulas for averages of characteristic polynomials the
function $z(x)$ plays the role for random matrix theory that
$\zeta(1+x)$ plays in the theory of moments of the Riemann
zeta-function.

Given finite subsets $A,B \subset \mathbb{C}$  we will have a sum
 over subsets $S\subset A$ and $T\subset B$ with $|S|=|T|$.  We let $\overline{S}=A-S$ and $\overline{T}=B-T$. We will
let $\hat \alpha$ denote a generic member of $S$ and $\hat \beta$
denote a generic member of $T$; we will use $\alpha$ and $\beta$ for
generic members of $A$ and $B$ or of $\overline{S}$ and
$\overline{T}$, according to the context. Also $S^-=\{-\hat\alpha:
\hat\alpha\in S\}$, and similarly for $T^-$. We let
\begin{eqnarray}
Z(A,B):=\prod_{\alpha\in A\atop\beta\in B}z(\alpha+\beta).
\end{eqnarray}
 
A simple modification of Theorem 4 of [CoSn] 
yields
\begin{theorem} \label{th:cosn}
Let $\delta>0$ and  let $\int_{(c)}$ denote an integration along the vertical path from $c-i\infty$ to $c+i\infty$. 
 Suppose that $F(x_1,\dots,x_n)$ is a holomorphic function which decays rapidly in each variable in horizontal strips.  Then, for any $\delta>0$,
\begin{eqnarray}\label{eq:offtheline}&&\int_{U(N)}\sum_{j_1,\dots ,j_n}
F(\theta_{j_1},\dots,\theta_{j_n})dX\nonumber
\\
&&\qquad =\frac{1}{(2\pi i)^n} \sum_{K+L+M=
\{1,\dots,n\}}(-1)^{|L|+|M|} N^{|M|} \\
&&\qquad\qquad\qquad \times\int_{(\delta)^{|K|}} \int_{(-\delta)^{|L|}}\int_{(0)^{|M|}}J^*(z_K;-z_L) F(iz_1,\dots,iz_n)~dz_1\dots
~dz_n\nonumber
\end{eqnarray}
where
  $z_K=\{z_k:k\in K\}$,  $-z_L=\{-z_\ell:\ell\in L\}$ and $\int_{(\delta)^{|K|}} \int_{(-\delta)^{|L|}}\int_{(0)^{|M|}}$
  means that we are integrating all the variables in $z_K$ along the $(\delta)$
  path,   all of the variables in $z_{L}$  along the $(-\delta)$
  path and all of the variables in $z_{M}$  along the $(0)$ path; and 
\begin{eqnarray} \label{eqn:Jsum} &&J^*(A,B):= \nonumber \\
 &&\qquad\qquad\sum_{S\subset A,T\subset B\atop
|S|=|T|}e^{-N(\sum_{\hat \alpha\in S} \hat \alpha
+\sum_{\hat{\beta}\in T}\hat\beta)} \frac{Z(S,T)Z(S^-,T^-)} {
Z^{\dagger}(S,S^-)Z^{\dagger}(T,T^-)} \sum_{{(A-S)+ (B-T)\atop =
U_1+\dots + U_Y}\atop |U_y|\le 2}\prod_{y=1}^Y H_{S,T}(U_y),
\end{eqnarray}
where
\begin{equation}\label{eqn:Hrmt}
H_{S,T}(W)=\left\{\begin{array}{ll} \sum_{\hat \alpha\in
S}\frac{z'}{z}(\alpha-\hat{\alpha})-\sum_{\hat\beta\in T}
\frac{z'}{z}(\alpha +\hat \beta) &\mbox{ if $W=\{\alpha\}\subset
A-S$}
   \\
\sum_{\hat\beta\in T}\frac{z'}{z}(\beta-\hat
\beta)-\sum_{\hat\alpha\in S} \frac{z'}{z}
(\beta+\hat\alpha) &\mbox{ if  $W=\{\beta\}\subset B-T$}\\
\left(\frac{z'}{z}\right)'(\alpha+\beta) & \mbox{ if
$W=\{\alpha,\beta\}$ with $
{\alpha \in A-S, \atop \beta\in B-T}$}\\
0&\mbox{ otherwise}.
\end{array}\right.
\end{equation}
The innermost sum of (\ref{eqn:Jsum}) is a sum over all partitions of $(A-S)+(B-T)$ into singletons and doubletons
$U_1, \cdots ,U_Y$. 
Also, $Z(A,B)=\prod_{\alpha\in A\atop\beta\in B}z(\alpha+\beta)$,
with the dagger imposing the additional restriction that a factor
$z(x)$ is omitted if its argument is zero:
\begin{eqnarray}
Z^\dagger(A,B)=\prod_{\alpha\in A\atop{\beta\in B\atop \alpha+\beta\ne 0}} z(\alpha+\beta).
\end{eqnarray}
  \end{theorem}
We remark that, assuming the  ratios conjecture, see [CoSn],  there is a structurally identical formula 
giving the $n$-correlation of  the zeros of the Riemann zeta-function.

A brief description of the proof of Theorem \ref{th:cosn} is that we start, in [CoSn], with a
theorem for the average of products of ratios of characteristic polynomials,
with shifts. Then one differentiates with respect to the shifts in the numerator to
get
a formula for the average of products of logarithmic derivatives of characteristic
polynomials. One then uses the latter formula in a residue computation to
obtain $n$-correlation sums.

  We want to apply the above theorem but with 
\begin{eqnarray} \label{eqn:testF}
F(x_1,\dots, x_n)=f\left(\frac{Nx_1}{2\pi},\dots,\frac{Nx_n}{2\pi}\right)h_1\left(\frac{x_1}{\mathcal{T}}\right)\dots h_n\left(\frac{x_n}{\mathcal{T}}\right)
\end{eqnarray}
 with $f\in \mathcal T_q$. 
 
In the right-hand side of the formula we replace $f$ by its Fourier transform so that we can better see what the 
implications of limited support are. Thus, we write
\begin{eqnarray*}
 f\left(\frac{iNz_1}{2\pi},\dots,\frac{iNz_n}{2\pi}\right)
=\int_{{\mathbf R}^n} \Phi(\xi_1,\dots,\xi_n)\delta(\xi_1+\dots+\xi_n) e\left(-\frac{iNz_1\xi_1 }{2\pi}-\dots-\frac{iNz_n\xi_n}{2\pi}\right)
d\xi_1\dots d\xi_n .
\end{eqnarray*}
Observe that 
\begin{eqnarray}
 e\left(-\frac{iNz_1\xi_1 }{2\pi}-\dots-\frac{iNz_n\xi_n}{2\pi}\right)
=e^{Nz_1\xi_1+\dots +Nz_n\xi_n}.
\end{eqnarray}
The contour integrals in Theorem \ref{th:cosn} restrict us to $|\Re z_i| \leq \delta$ for each $i$.
  Suppose that $\Phi(\xi_1,\dots ,\xi_n)=0$ if $|\xi_1|+\dots+ |\xi_n|>2q-\epsilon$ for some $\epsilon>0$. Then $|\Phi|\ne0$ implies that 
\begin{eqnarray}
\left|e^{Nz_1\xi_1+\dots +Nz_n\xi_n}\right|\le e^{N\delta(2q-\epsilon)}.
\end{eqnarray}

 We compare this exponential with the exponentials which appears in the factor  $J^*(z_K;-z_L)$:  
\begin{eqnarray}
\sum_{S\subset z_K,T\subset -z_L\atop
|S|=|T|}e^{-N(\sum_{\hat \alpha\in S} \hat \alpha
+\sum_{\hat{\beta}\in T}\hat\beta)}.
\end{eqnarray}
Notice that the real parts of all of the $\alpha\in z_K$ and all of the $\beta\in -z_L$ are equal to $\delta$.  If $|S|=|T|\ge q$, then
\begin{eqnarray}
\left|e^{-N(\sum_{\hat \alpha\in S} \hat \alpha
+\sum_{\hat{\beta}\in T}\hat\beta)}\right|\le e^{-2N\delta q}.
\end{eqnarray}
Thus, the product of these two factors is $\le e^{-N\delta \epsilon}$.
We can move the paths of integration $(-\delta)$ and $(\delta)$ away from the imaginary axis. 
 The integrand tends to zero uniformly on the vertical paths as $\delta \rightarrow \infty$. 
Note that easy estimates show that the other factors in the integrand do not 
interfere with this estimate. For example, if $\alpha\in A, \beta\in B$ then $\Re(\alpha+\beta)=2\delta$ so that
\begin{eqnarray}
|z(\alpha+\beta)|\le (e^{2\delta}-1)^{-1} \to 0
\end{eqnarray}
as $\delta\to \infty$. Also, if $\hat{\alpha}\in S$ and $\alpha\in S^-$ then $\Re (\alpha+\hat{\alpha})=0$ so that 
\begin{eqnarray}\frac{1}{|z(\alpha+\hat{\alpha})|}\le 2.
\end{eqnarray}   
The $H$-terms involve logarithmic derivatives of $z$ for which the identity $z'/z=1-z$ is helpful.
Finally the $h$ functions can be bounded by
\begin{eqnarray}
h(iz/\mathcal{T})|\ll e^{\Delta \delta/\mathcal{T}}
\end{eqnarray}
where $\Delta$ is such that $g$ is supported on $[-\Delta,\Delta]$. Thus, the product of the $h$ may be bounded by $e^{n\Delta \delta/\mathcal{T}}$.
Since $N$ is assumed sufficiently large with respect to $n$ (and $\mathcal{T} \to \infty$ before $N$) we see that the 
terms with $|S|=|T|\geq q$ are 0.

We let $J_q^*(A;B)$ be defined as $J^*(A;B)$ but with the subsets $S$ and $T$ in the defining sum having size smaller than $q$, i.e.
\begin{eqnarray}
J_q^*(A;B)= \sum_{S\subset A, T\subset B\atop |S|=|T| < q} \dots .
\end{eqnarray}
Then, if the total support of $\Phi$ is limited to any number smaller than $2q$,  then Theorem \ref{th:cosn}
 holds with all the $J^*$ replaced by $J_q^*$.

\begin{theorem} \label{th:restrictedF}
Let $\delta>0$ and  suppose that $F(x_1,\dots,x_n)$ satisfies (\ref{eqn:testF}).  Then,  
\begin{eqnarray} &&\int_{U(N)}\sum_{j_1,\dots ,j_n}
F(\theta_{j_1},\dots,\theta_{j_n})dX\nonumber
\\
&&\qquad =\frac{1}{(2\pi i)^n} \sum_{K+L+M=
\{1,\dots,n\} }(-1)^{|L|+|M| } N^{|M|} \\
&&\qquad\qquad\qquad \times\int_{(\delta)^{|K|}} \int_{(-\delta)^{|L|}}\int_{(0)^{|M|}}J_q^*(z_{K};-z_{L}) F(iz_1,\dots,iz_n)~dz_1\dots
~dz_n\nonumber
\end{eqnarray}
\end{theorem}

\subsection{Examples.} We give some examples of $J^*$ and $J^*_q$ to help the reader parse the 
last two theorems.  We write
\begin{eqnarray}
J^*(A,B)=\sum_{{S\subset A\atop T\subset B}\atop |S|=|T|}
D_{S,T} (A-S,B-T)
\end{eqnarray}
with an obvious notation. We have
\begin{eqnarray}
&&J^*(\{a\},\{b\})= 
D_{\phi,\phi}(\{a\},\{b\})+D_{\{a\},\{b\}}(\phi,\phi),
\end{eqnarray}
where
\begin{eqnarray}
D_{\phi,\phi}(\{a\},\{b\})=\left(\frac{z'}{z}\right)'(a+b),
\end{eqnarray}
and
\begin{eqnarray}
D_{\{a\},\{b\}}(\phi,\phi)=e^{-N(a+b)}z(a+b)z(-a-b).
\end{eqnarray}
Thus,
\begin{eqnarray}
J_1(\{a\},\{b\})=\left(\frac{z'}{z}\right)'(a+b)
\end{eqnarray}
and if $q\geq 2$, then 
\begin{eqnarray}
J_q(\{a\},\{b\})=J(\{a\},\{b\})=e^{-N(a+b)}z(a+b)z(-a-b)+\left(\frac{z'}{z}\right)'(a+b).
\end{eqnarray}
Next,
\begin{eqnarray}
&&J^*(\{a\},\{b_1,b_2\})=  D_{\phi,\phi}(\{a\},\{b_1,b_2\})+ D_{\{a\},\{b_1\}}(\phi,\{b_2\})+D_{\{a\},\{b_2\}}(\phi,\{b_1\}),
\end{eqnarray}
where
\begin{eqnarray}
D_{\phi,\phi}(\{a\},\{b_1,b_2\})=0,
\end{eqnarray}
\begin{eqnarray}
D_{\{a\},\{b_1\}}(\phi,\{b_2\})=e^{-N(a+b_1)}z(a+b_1)z(-a-b_1)
\left(\frac{z'}{z}(b_2-b_1)-\frac{z'}{z}(b_2+a) \right),
\end{eqnarray}
and
\begin{eqnarray}
D_{\{a\},\{b_2\}}(\phi,\{b_1\})=e^{-N(a+b_2)}z(a+b_2)z(-a-b_2)
\left(\frac{z'}{z}(b_1-b_2)-\frac{z'}{z}(b_1+a) \right),
\end{eqnarray}
Thus, $J_1^*(\{a\},\{b_1,b_2\})=0$ and if $q\geq 2$, then
\begin{eqnarray*} 
J_q^*(\{a\},\{b_1,b_2\})
  &=& J^*(\{a\},\{b_1,b_2\})\\
&=&
e^{-N(a+b_1)}z(a+b_1)z(-a-b_1)
\left(\frac{z'}{z}(b_2-b_1)-\frac{z'}{z}(b_2+a) \right)\\
&& \qquad +  \nonumber e^{-N(a+b_2)}z(a+b_2)z(-a-b_2)
\left(\frac{z'}{z}(b_1-b_2)-\frac{z'}{z}(b_1+a) \right) .
\end{eqnarray*}
Next
\begin{eqnarray*}  \label{eqn:Jabbb}
&& J^*(\{a\},\{b_1,b_2,b_3\})=D_{\phi,\phi}(\{a\},\{b_1,b_2,b_3\})+
D_{\{a\},\{b_1\}}(\phi,\{b_2,b_3\})+D_{\{a\},\{b_2\}}(\phi,\{b_1,b_3\})\nonumber \\
&& \qquad \qquad \qquad \qquad \qquad +D_{\{a\},\{b_3\}}(\phi,\{b_1,b_2\})\nonumber
\\
&&  \qquad = e^{-N(a+b_1)}z(a+b_1)z(-a-b_1)
\left(\frac{z'}{z}(b_2-b_1)-\frac{z'}{z}(b_2+a) \right)
\left(\frac{z'}{z}(b_3-b_1)-\frac{z'}{z}(b_3+a) \right)  \nonumber\\
&& \qquad \quad +    e^{-N(a+b_2)}z(a+b_2)z(-a-b_2)
\left(\frac{z'}{z}(b_1-b_2)-\frac{z'}{z}(b_1+a)
   \right) \left(\frac{z'}{z}(b_3-b_2)-\frac{z'}{z}(b_3+a) \right)
\\
&& \qquad \quad +   \nonumber e^{-N(a+b_3)}z(a+b_3)z(-a-b_3)
\left(\frac{z'}{z}(b_1-b_3)-\frac{z'}{z}(b_1+a)
   \right) \left(\frac{z'}{z}(b_2-b_3)-\frac{z'}{z}(b_2+a) \right)
\end{eqnarray*}

With  $A=\{a_1, a_2\},B=\{b_1,b_2\}$  we have
\begin{eqnarray*}
&&J^*(\{a_1,a_2\},\{b_1,b_2\})=  D_{\phi,\phi}(\{a_1, a_2\},\{b_1,b_2\})+ D_{\{a_1\},\{b_1\}}(\{a_1\},\{b_1\})\\
&& \qquad \qquad \qquad \qquad \qquad+
D_{\{a_1\},\{b_2\}}(\{a_2\},\{b_1\})+
D_{\{a_2\},\{b_1\}}(\{a_1\},\{b_2\})\\
&& \qquad \qquad \qquad \qquad \qquad +D_{\{a_2\},\{b_2\}}(\{a_1\},\{b_1\})+D_{\{a_1,a_2\},\{b_1,b_2\}}(\phi,\phi).
\end{eqnarray*}
Here the $J_1^*$ term will have only $D_{\phi,\phi}$, the $J_2^*$ term will have all but the $D_{\{a_1,a_2\},\{b_1,b_2\}}$
term and if $q\ge 3$ all of the above terms will be present in $J_q^*$, where
\begin{eqnarray*}
D_{\phi,\phi}(\{a_1,a_2\},\{b_1,b_2\})=\left(\frac{z'}{z}\right)'(a_1+b_1)\left(\frac{z'}{z}\right)'(a_2+b_2)
+\left(\frac{z'}{z}\right)'(a_1+b_2)\left(\frac{z'}{z}\right)'(a_2+b_1)
\end{eqnarray*}
and
\begin{eqnarray*}
D_{\{a_1\},\{b_1\}}(\{a_2\},\{b_2\}) &=&e^{-N(a_1+b_1)}z(a_1+b_1)z(-a_1-b_1)(H_{a_2,b_2}+H_{a_2}
H_{b_2})\\
&=& \nonumber e^{-N(a_1+b_1)}z(a_1+b_1)z(-a_1-b_1)
\left(\left(\frac{z'}{z}\right)'(a_2+b_2)\right.
\\
&&\qquad \nonumber\left.+\left(\frac{z'}{z}(a_2-a_1)-\frac{z'}{z}(a_2+b_1)
\right) \left(\frac{z'}{z}(b_2-b_1)-\frac{z'}{z}(b_2+a_1) \right)
\right);
\end{eqnarray*}
the other $D_{\{a_i\},\{b_j\}}$ are similar. Also,
\begin{eqnarray} &&
D_{\{a_1,a_2\},\{b_1,b_2\}}(\phi,\phi)=e^{-N(a_1+a_2+b_1+b_2)}\times \nonumber\\ \nonumber
&&\quad \frac{z(a_1+b_1)z(-a_1-b_1)z(a_1+b_2)z(-a_1-b_2)
z(a_2+b_1)z(-a_2-b_1)z(a_2+b_2)z(-a_2-b_2)}
{z(a_1-a_2)z(-a_2-a_1)z(b_1-b_2)z(b_2-b_1)}.
\end{eqnarray}
  
See the last section of [CoSn] for more detailed examples.

\section{The special case $q=1$}

If $q=1$, as in the theorem of Rudnick and Sarnak, then the sets $S$ and $T$ in the sum defining $J_1^*(A;B)$ are both empty. We have
\begin{eqnarray}
J_1^*(A;B)=
 \sum_{{A+B =\atop 
U_1+\dots + U_Y}\atop |U_y|\le 2}\prod_{y=1}^Y H_{\emptyset,\emptyset}(U_y),
\end{eqnarray}
We have
\begin{equation}\label{eqn:Hrmtrestricted}
H_{\emptyset,\emptyset}(U)=\left\{\begin{array}{ll}  
\left(\frac{z'}{z}\right)'(\alpha+\beta) & \mbox{ if
$U=\{\alpha,\beta\}$ with $
{\alpha \in A, \atop \beta\in B}$}\\
0&\mbox{ otherwise}.
\end{array}\right.
\end{equation}
In particular, in the sum $A+B=U_1+\dots+U_Y$ we have that $|U_y|=2$ for each $y$ with each $U_y$ having precisely 
one element from $A$ and one element from $B$. In particular, $|A|=|B|$.
This gives, when $A=\{\alpha_1,\dots,\alpha_{|K|}\},B=\{\beta_1,\dots,\beta_{|K|}\}$, 
\begin{equation}
J_1^*(A;B)=\sum_{\sigma\in S_{|K|}}\prod_{k=1}^{|K|} 
\left(\frac{z'}{z}\right)'(\alpha_k+\beta_{\sigma(k)}),
\end{equation}
and 0 otherwise. 
  Thus, we have 
\begin{theorem} \label{th:q1}
Let $\delta>0$ and  suppose that $F(x_1,\dots,x_n)$ satisfies (\ref{eqn:testF}) with $q=1$.  Then,  
\begin{eqnarray} &&\int_{U(N)}\sum_{j_1,\dots ,j_n}
F(\theta_{j_1},\dots,\theta_{j_n})dX\nonumber
\\
&&\qquad =\frac{1}{(2\pi i)^n} \sum_{K+L+M=
\{1,\dots,n\}\atop |K|=|L|}(-1)^{|L|+|M|} N^{|M|} \\
&&\qquad\qquad\qquad \times\int_{(\delta)^{|K|}} \int_{(-\delta)^{|L|}}\int_{(0)^{|M|}}\sum_{\sigma\in S_{|K|}}\prod_{j=1}^{|K|}
\left(\frac{z'}{z}\right)'(z_{k_j}-z_{\ell_{\sigma(j)}})F(iz_1,\dots,iz_n)~dz \nonumber
\end{eqnarray}
where 
\begin{eqnarray} K=\{k_1,\dots,k_{|K|}\} \qquad \mbox{and} \qquad L=\{\ell_1,\dots,\ell_{|K|}\}.
\end{eqnarray}
\end{theorem}

\section{Comparison with Rudnick-Sarnak}

We now put the formula from the last section into the form of Theorem 3.1 of [RuSa]. For ease of subscripting, consider one particular term, 
denoted (with an abuse of notation that hopefully won't cause confusion)
$K=\{1,2,\dots, K\}$, $L=\{K+1,\dots ,2K\}$ and $\sigma$ is the identity permutation. Let 
\begin{eqnarray}
I:= \int_{(\delta)^{K}}\int_{(-\delta)^{K}}  \int_{(0)^{n-2K}} \prod_{k=1}^K 
\left(\frac{z'}{z}\right)'(z_k-z_{K+k})F(iz_1,\dots,iz_n)~dz_1\dots
~dz_n.\nonumber
\end{eqnarray}
 We write 
\begin{eqnarray}&&
F(iz_1,\dots,iz_n) =h_1\left(\frac{iz_1}{\mathcal{T}}\right)\dots h_n\left(\frac{iz_n}{\mathcal{T}}\right)\int_{\mathbf R^n}\Phi(\xi_1,\dots,\xi_n)\\
&&\qquad 
\times \delta(\xi_1+\dots +\xi_n) e\left(-\frac{iNz_1\xi_1 }{2\pi }-\dots -\frac{iNz_n\xi_n}{2\pi }\right)
~d\xi_1\dots d\xi_n.\nonumber
\end{eqnarray}
We integrate the $z_r$ variables with $2K<r\le n$ using equation (\ref{eqn:gtoh}) and Fourier inversion to get
\begin{eqnarray}&&\nonumber
\int_{(0)^{n-2K}}
e\left(-\frac{iNz_{2K+1}\xi_{2K+1} }{2\pi }-\dots -\frac{iNz_n\xi_n}{2\pi }\right)
h_{2K+1}\left(\frac{iz_{2K+1}}{\mathcal{T}}\right)\dots h_n\left(\frac{iz_n}{\mathcal{T}}\right)~dz_{2K+1}\dots dz_n\\
&&\qquad =\prod_{j=2K+1}^n \int_{(0)}h_j\left( \frac{iz_j}{\mathcal{T}}\right) e\left(-\frac{iNz_j\xi_j}{2\pi }\right) dz_j
=\prod_{j=2K+1}^n\left(-2\pi i\mathcal{T} g_j(N\mathcal{T}\xi_j)\right).
\end{eqnarray}
This gives 
\begin{eqnarray}&& \nonumber
I=\int_{\mathbf R^n}\Phi(\xi_1,\dots,\xi_n)\delta(\xi_1+\dots +\xi_n)\prod_{j=2K+1}^n\left(-2\pi i\mathcal{T} g_j(N\mathcal{T}\xi_j)\right)
 \int_{(\delta)^{K}}  \int_{(-\delta)^{K}}  \prod_{k=1}^{K}
\\ && \qquad \left(
 e\left(\frac{-iNz_k\xi_k}{2\pi }\right)h_k\left(\frac{iz_k}{\mathcal{T}}\right) e\left(-\frac{iNz_{k+K}\xi_{K+k}}{2\pi }\right)\right.\\
&&\qquad\qquad\qquad \times \left.
h_{K+k}\left(\frac{iz_{K+k}}{\mathcal{T}}\right) 
\left(\frac{z'}{z}\right)'(z_k-z_{K+k})\right) ~dz_1\cdots dz_{2K} d\xi.
\nonumber  \end{eqnarray}
Now we move the paths of integration of the $z_{k}$ with $1\le k\le K$. If $\xi_{k}>0$ we move the path to the left 
across the double pole at $z_{K+k}$; if $\xi_{k}<0$ then we move the path far to the right. We use the fact that 
\begin{eqnarray}
\operatornamewithlimits{Res}_{z_k=z_{K+k}}\frac{f_1(z_k)f_2(z_{K+k})}{(z_k-z_{K+k})^2}=f_1'(z_{K+k}) f_2(z_{K+k}).
\end{eqnarray}
Setting
$ f_1(z_k)=e(-i N z_k \xi_k/2 \pi) h_k(iz_k/\mathcal{T})$ and $ f_2(z_{k+K})=e(-i N z_{k+K} \xi_{k+K}/2 \pi) h_{k+K}(iz_{k+K}/\mathcal{T})$
and applying the above identity we get 
\begin{eqnarray}&&\nonumber
I=\left(1+O(1/\mathcal{T})\right)(2\pi i)^K\int_{\mathbf R^n\atop \xi_j>0, j\le K}\Phi(\xi_1,\dots,\xi_n)\delta(\xi_1+\dots +\xi_n)\prod_{j=2K+1}^n
\left(-2\pi i\mathcal{T}g_j\left(  N\mathcal{T}\xi_j\right)\right)
 \\ && \quad \times
\int_{(-\delta)^{K}}  \prod_{k=1}^K\left(N\xi_{k}
 e\left(\frac{-iNz_{K+k}(\xi_k+\xi_{K+k})}{2\pi }\right)h_k\left(\frac{iz_{K+k}}{\mathcal{T}}\right) 
h_{K+k}\left(\frac{iz_{{K+k}}}{\mathcal{T}}\right) 
  ~dz_{K+k}\right) d\xi.
\end{eqnarray}
We compute, for example,  
\begin{eqnarray*}
\int_{(-\delta)}  
 e\left(-\frac{iNz(\xi_1+\xi_{2})}{2\pi }\right)h_1\left(\frac{iz}{\mathcal{T}}\right) 
h_{2}\left(\frac{iz}{\mathcal{T}}\right) 
  ~dz =-2\pi i\mathcal{T}\int_{\mathbf R} g_1(u)g_2(-u+N\mathcal{T}(\xi_1+\xi_2))~du
\end{eqnarray*}
Thus, 
\begin{eqnarray}&&\nonumber
I=(2\pi iN)^K\int_{ \mathbf R^n\atop \xi_j>0, j\le K }
\xi_{1}\dots \xi_{K}\Phi(\xi_1,\dots,\xi_n)\delta(\xi_1+\dots +\xi_n)\prod_{j=2K+1}^n\left(-2\pi i\mathcal{T}g_j(N\mathcal{T}\xi_j) \right)
 \\ && \qquad \times
  \prod_{k=1}^K \left(-2\pi i\mathcal{T} 
 \int_{\mathbf R}g_k(u_k)g_{K+k}(-u_k+\mathcal{T}N(\xi_k+\xi_{K+k})) ~du_k\right)
  d\xi(1+O(\tfrac{1}{\mathcal{T}})).
\end{eqnarray}
We make the changes of variables 
\begin{eqnarray}
\begin{array}{lll}
y_j&=N\mathcal{T}\xi_j   &\mbox{for $2K+1\le j\le n$}\\
y_k&=u_k  &\mbox{for $1\le k\le K$}\\
y_{K+k} &= -u_k + N\mathcal{T}(\xi_k+\xi_{K+k})   &\mbox{for $1\le k\le K$}.
\end{array}
\end{eqnarray} 
The last substitution implies that 
\begin{eqnarray}
\xi_{K+k}=-\xi_{k}+\frac{ (y_{K+k}+u_k) }{N\mathcal{T}}.
\end{eqnarray}
Also, the condition $\sum_{j=1}^n \xi_j=0$ implied by the delta-function translates to $\sum_{j=1}^n y_j=0$. We have 
\begin{eqnarray}&&\nonumber
I = \frac{(2\pi iN)^{K} ( -2\pi i \mathcal{T})^{n-K}}{(N\mathcal{T})^{n-K-1}}\int_{\mathbf R^n\atop \sum_{y_j=0}}\prod_{j=1}^n g_j(y_j) \int_{  \xi_{k}>0, 1\le k\le K}
\xi_{1}\dots \xi_{K}
\\ &&\qquad \times
\Phi(\xi_1,\dots, \xi_k,-\xi_1+\frac{y_{K+1}+y_1}{N\mathcal{T}},\ldots,-\xi_{K}+\frac{ (y_{2K}+y_K) }{N\mathcal{T}},\frac{y_{n-2K}}{N\mathcal{T}},\dots, \frac{ y_{n}}{N\mathcal{T}})\\
&& \qquad \qquad \times (1+O(\tfrac{1}{\mathcal{T}}))
d\xi_{1}\dots d\xi_{K} ~dy_1\cdots dy_n. \nonumber
\end{eqnarray}
Employing the Taylor expansion of $\Phi$, just as in [RuSa], we have
\begin{eqnarray} && \nonumber
\frac{I}{N\mathcal{T}}=N^{2K-n} ( 2\pi i )^{n}(-1)^{n-K}
\int_{\mathbf R^n\atop \sum_{y_j=0}}\prod_{j=1}^n g_j(y_j)~dy \int_{\mathbf R^K\atop \xi_{k}>0, 1\le k\le K}
\xi_{1}\dots \xi_{K}
\\ &&\qquad \times
\Phi(\xi_{1},\dots,\xi_{K},-\xi_{1},\dots,-\xi_{K}, 0,\dots,0)  ~d\xi_{1}\dots d\xi_{K}(1+O(1/\mathcal{T})).
\end{eqnarray}
Since 
\begin{eqnarray}\nonumber
\int_{\mathbf R^n\atop \sum_{y_j=0}}\prod_{j=1}^n g_j(y_j)~dy&=&\int_{\mathbf R^{n-1}}
\prod_{j=1}^{n-1} g_j(y_j)\frac{1}{2\pi}\int_{\mathbf R} h_n(t)\exp(-it(\sum_{j=1}^{n-1} y_j) ~dt ~dy\\
&=&\frac{1}{2\pi}\int_{\mathbf R} \prod_{j=1}^n h_j(t)~dt =
\frac{\kappa(\mathbf h)}{2\pi},
\end{eqnarray}
 we have
\begin{eqnarray} &&\nonumber
\frac{I}{N\mathcal{T}}= N^{2K-n} ( 2\pi i )^{n}(-1)^{n-K}
\frac{\kappa(\mathbf h)}{2\pi}\int_{\mathbf R^K\atop \xi_{k}>0, 1\le k\le K}
\xi_{1}\dots \xi_{K}
\\ &&\qquad \times
\Phi(\xi_{1},\dots,\xi_{K},-\xi_{1},\dots,-\xi_{K}, 0,\dots,0)  ~d\xi_{1}\dots d\xi_{K}(1+O(1/\mathcal{T})).
\end{eqnarray}
More generally, for $K=\{k_1,\dots,k_{|K|}\}, L=\{\ell_1,\dots,\ell_{|K|}\}$ and $\sigma\in S_{|K|}$,
\begin{eqnarray}
I(K,L,\sigma):&=&\int_{(\delta)^{|K|}} \int_{(-\delta)^{|L|}}\int_{(0)^{|M|}}  \prod_{j=1}^K
\left(\frac{z'}{z}\right)'(z_{k_j}-z_{\ell_{\sigma(j)}})F(iz_1,\dots,iz_n)~dz \nonumber\\
&\sim&(N\mathcal{T})N^{2|K|-n} ( 2\pi i )^{n}(-1)^{n-|K|}
\frac{\kappa(\mathbf h)}{2\pi}\int_{  \xi_{k_j}>0 }
\xi_{k_1}\dots \xi_{k_{|K|}}
 \Phi\left(\sum_{j=1}^{|K|} \xi_{k_j}\mathbf e_{k_j,\ell_{\sigma(j)}}\right) d\xi .
\end{eqnarray}
We insert this into Theorem \ref{th:q1} and have
\begin{eqnarray} &&\int_{U(N)}\sum_{j_1,\dots ,j_n}
F(\theta_{j_1},\dots,\theta_{j_n})dX\nonumber
\\
&&\qquad =\frac{1}{(2\pi i)^n} \sum_{K+L+M=
\{1,\dots,n\}\atop |K|=|L|}(-1)^{|L|+|M|} N^{|M|} N^{2|K|-n} ( 2\pi i )^{n}(-1)^{n-|K|}\\
&&\qquad\qquad\qquad \times 
\frac{\kappa(\mathbf h)N\mathcal{T}}{2\pi}\int_{  \xi_{k_j}>0 }
\xi_{k_1}\dots \xi_{k_{|K|}}\sum_{\sigma \in S_K}
 \Phi\left(\sum_{j=1}^{|K|} \xi_{k_j}\mathbf e_{k_j,\ell_{\sigma(j)}}\right) d\xi(1+O(1/\mathcal{T})).\nonumber
\end{eqnarray}
Since $|M|=n-2|K|$ the right-hand side simplifies to 
\begin{eqnarray} 
\kappa(\mathbf h)\frac{N\mathcal{T}}{2\pi}  \sum_{K+L+M=
\{1,\dots,n\}\atop |K|=|L|} \sum_{\sigma\in S_{|K|}} 
\int_{  \xi_{k_j}>0 }
\xi_{k_1}\dots \xi_{k_{|K|}}
 \Phi\left(\sum_{j=1}^{|K|} \xi_{k_j}\mathbf e_{k_j,\ell_{\sigma(j)}}\right) d\xi +O(N).
\end{eqnarray}
This proves Theorem \ref{th:r1}.

\section{Acknowledgment}
The authors would like to thank the referee for helpful suggestions.

\end{document}